%%%%%%%%%%%%%%%%%%%%%%%%%%%%%%%%%%%%%%%%%%%%%%%%%%%%%%%%%%%%%%%%%%%%%%%%%%%
%%% Title:  Pseudoriemannian metrics on spaces of bilinear structures 
%%% Author: Olga Gil-Medrano, Peter W. Michor, Martin Neuwirther
%%% Remark: AmSTeX, 21 pages
%%% Series: Quarterly J\. Math Oxford (2), 43 (1992), 201--221
%%%%%%%%%%%%%%%%%%%%%%%%%%%%%%%%%%%%%%%%%%%%%%%%%%%%%%%%%%%%%%%%%%%%%%%%%%%
\input amstex
\input amsppt.sty
\def\nmb#1#2{#2}         % used for renumbering, TeX should ignore.
\def\totoc{}             %= to table of content, invoked by amspptb.sty
\def\idx{}               % for producing index, invoked by amspptb.sty
\def\ign#1{}             %=ignore, invisible entry for the index only
\pagewidth{11.8cm} 
\pageheight{18cm}
\magnification=\magstep1

%\pageno201
\redefine\o{\circ}
\define\X{\frak X}
\define\al{\alpha}

\define\ga{\gamma}
\define\de{\delta}

\define\ze{\zeta}
\define\et{\eta}

\define\ka{\kappa}
\define\la{\lambda}

\define\ph{\varphi}

\define\ps{\psi}
\define\om{\omega}
\define\Ga{\Gamma}

\define\La{\Lambda}

\define\Ph{\Phi}
\define\Ps{\Psi}
\define\Om{\Omega}
\predefine\ii{\i}
\redefine\i{^{-1}}
\define\row#1#2#3{#1_{#2},\ldots,#1_{#3}}
\define\x{\times}
\define\M{\Cal M}
\redefine\D{\Cal D}
\define\Vol{\operatorname{vol}}
\define\Diff{\operatorname{Diff}}
\define\tr{\operatorname{tr}}
\define\dd#1{\frac{\partial}{\partial#1}}
\define\Exp{\operatorname{Exp}}

\define\artanh{\operatorname{Artanh}}
\define\im{\operatorname{im}}
\def\today{\ifcase\month\or
 January\or February\or March\or April\or May\or June\or
 July\or August\or September\or October\or November\or December\fi
 \space\number\day, \number\year}
%\newsymbol\lessgtr 1337
%\newsymbol\gtrless 133F
%\newsymbol\twoheadrightarrow 1310
\topmatter
\title Pseudoriemannian metrics on spaces \\
of bilinear structures \endtitle
\author  Olga Gil-Medrano\footnote{Research partially
supported the CICYT grant n. PS87-0115-G03-01, work done during a
stay at Vienna, September 1990 \hfill\hfill\hfill\hfill} \\
Peter W. Michor  \\
Martin Neuwirther\footnote{Supported by Project P 7724 PHY
of `Fonds zur F\"orderung der wissenschaftlichen Forschung'\hfill}
   \endauthor
%\affil
%Departamento de Geometr\'\ii a y Topolog\'\ii a\\
%Universidad de Valencia, Spain\\
%\\
%Institut f\"ur Mathematik, Universit\"at Wien, Austria.
%\endaffil
\leftheadtext{\smc Olga Gil-Medrano, Peter W. Michor, Martin Neuwirther}
\rightheadtext{\smc Metrics on spaces of bilinear structures}
\address{O\. Gil-Medrano:
Departamento de Geometr\'\ii a y Topolog\'\ii a,
Facultad de Matem\'a\-ti\-cas,
Universidad de Valencia,
46100 Burjassot,
Valencia, Spain}
\endaddress 
\address{P\. W\. Michor, M\. Neuwirther:
Institut f\"ur Mathematik, Universit\"at Wien,
Strudlhofgasse 4, A-1090 Wien, Austria}
\endaddress
 
%\date{\today}\enddate
\abstract{The space of all non degenerate bilinear structures on a
manifold $M$ carries a one parameter family of pseudo Riemannian
metrics. We determine the geodesic equation, covariant derivative,
curvature, and we solve the geodesic equation
explicitly. Each space of pseudo Riemannian metrics with
fixed signature is a geodesically closed submanifold. The space of non
degenerate 2-forms is also a geodesically closed submanifold.
Then we show that, if we fix a distribution on $M$, the space of all Riemannian
metrics splits as the product of three spaces which are
everywhere mutually orthogonal, for the usual metric. We
investigate this situation in detail.}
\endabstract
\endtopmatter
%\input amspptb.sty
%\userunningheads
%\def\leftheadtext{\smc Olga Gil-Medrano, Peter W. Michor, Martin Neuwirther}
%\def\rightheadtext{\smc Metrics on spaces of bilinear structures}
%\def\bottremark{\today\hfill}

\document
 
\heading Table of contents \endheading
%\input \jobname.toc
%\loadtoc
\noindent 0. Introduction \leaders \hbox to 1em{\hss .\hss }\hfill {\eightrm 1}\par 
\noindent 1. The general setup \leaders \hbox to 1em{\hss .\hss }\hfill {\eightrm 2}\par 
\noindent 2. Geodesics, Levi Civita connection, and curvature \leaders \hbox to 1em{\hss .\hss }\hfill {\eightrm 5}\par 
\noindent 3. Some submanifolds of $\Cal B$ \leaders \hbox to 1em{\hss .\hss }\hfill {\eightrm 11}\par 
\noindent 4. Splitting the manifold of metrics \leaders \hbox to 1em{\hss .\hss }\hfill {\eightrm 13}\par 

\heading\totoc 0. Introduction \endheading
 
If $M$ is a (not necessarily compact) smooth finite dimensional 
manifold, the space $\Cal B = C^\infty(GL(TM,T^*M))$ of all non 
degenerate $\binom02$-tensor fields on it can be endowed with a 
structure of an infinite dimensional smooth manifold modeled on the 
space $C_c^\infty(L(TM,T^*M))$ of $\binom02$ -tensor fields with 
compact support, in the sense of \cite{Michor, 1980}. The tangent 
bundle of $\Cal B$ is $T\Cal B = \Cal B\x\ C_c^\infty(L(TM,T^*M))$ 
and we consider on $\Cal B$ the one parameter family of smooth 
pseudo Riemannian metrics $G^\al$ defined in section 1.   

Section 2 is devoted to the study of the geometry of $(\Cal B,G^\al)$. 
We start by computing the geodesic equation, using variational 
methods, and from it the covariant derivative and the curvature of 
the Levi-Civita connection are obtained. We explicitly solve the 
geodesic equation and we find the domain of definition of the 
exponential mapping which is open for the topology considered on 
$\Cal B$. We show that the exponential mapping is a real analytic 
diffeomorphism from an open neighborhood of zero onto its image.   

Pseudo Riemannian metrics of fixed signature and non degenerate 
2-forms are splitting submanifolds of $\Cal B$. In section 3 we show 
that they are geodesically closed. In particular this applies to the 
manifold $\M$ of all Riemannian metrics on $M$, with its usual 
metric, the geometry of which has been studied in \cite{Ebin, 1970}, 
\cite{Freed, Groisser, 1989}, \cite{Gil-Medrano, Michor, 1990}. If $M$ 
is compact, the space of symplectic structures on $M$ is also a 
splitting submanifold, it seems that it should admit a connection for 
the induced metric, but we have not been able to find it because the 
exterior derivative has a complicated expression in the terms in 
which the metric is simple.   

Finally in section 4 we investigate the 
splitting of $\M$ induced by fixing a distribution $V\subset TM$: it 
turns out to be the product of the space $\M(V)$ of all fiber metrics 
on the distribution, the space $\M(TM/V)$ of all fiber metrics on the 
normal, and the space $\Cal P_V(M)$ of all almost product structures 
having $V$ as the vertical distribution. Each slice 
$\M(V)\x \M(TM/V)$ is a geodesically closed submanifold and the 
induced metric on each slice $\Cal P_V(M)$ is flat. As a Riemannian 
manifold $\Cal M$ turns out to be what is usually known as a product 
with varying metric on the fibers. If the distribution $V$ is the 
vertical bundle of a fiber bundle $(E,p,B)$ then $\Cal P_V(E)$ is 
just the space of all connections and we expect that results of this 
section could be used for studying the moduli space of connections 
modulo the gauge group.  This paper is a sequel to 
\cite{Gil-Medrano, Michor, 1990} and essentially the same techniques 
and ideas are used. We will refer to it frequently. 
 
\heading \totoc\nmb0{1}. The general setup \endheading
 
\subheading{\nmb.{1.1}. Bilinear structures}
Let $M$ be a smooth second countable
finite dimensional manifold. Let $\otimes^2T^*M$ denote the
vector bundle of all
$\binom02$-tensors on $M$
which we canonically identify with the bundle $L(TM,T^*M)$. Let
$GL(TM,T^*M)$ denote the non degenerate ones.
For any $b:T_xM\to T_x^*M$ we let the \idx{\it transposed} be given
by $b^t:T_xM\to T_x^{**}M @>b'>> T_x^*M$. As a bilinear structure $b$
is skew symmetric if and only if $b^t=-b$, and $b$ is symmetric if
and only if $b^t=b$. In the latter case a frame $(e_j)$ of $T_xM$ can be
chosen in such a way that in the dual frame $(e^j)$ of $T^*_xM$ we
have
$$b = e^1\otimes e^1 + \dots + e^p\otimes e^p
     - e^{p+1}\otimes e^{p+1}-e^{p+q}\otimes e^{p+q};$$
$b$ has signature $(p,q)$ and is non degenerate if and only if
$p+q=n$, the dimension of $M$. In this case $q$ alone will be called 
the signature.
 
A section $b\in C^\infty(GL(TM,T^*M))$ will be called a non
degenerate bilinear structure on $M$ and we will denote the space of
all such structures by
$\Cal B(M)=\Cal B:=C^\infty(GL(TM,T^*M))$.
It is open in the
space of sections $C^\infty(L(TM,T^*M))$ for the Whitney
$C^\infty$-topology, in which the latter space is, however, not
a topological vector space, since $\frac1n h$ converges to 0 if
and only if $h$ has compact support. So the space
$\Cal B_c=C^\infty_c(L(TM,T^*M))$
of sections with compact support is the largest topological
vector space contained in the topological group
$(C^\infty(L(TM,T^*M)),+)$, and the trace of the Whitney
$C^\infty$-topology on it coincides with the inductive limit topology
$$C^\infty_c(L(TM,T^*M))=\varinjlim_KC^\infty_K(L(TM,T^*M)),$$
where $C^\infty_K(L(TM,T^*M))$ is the space of all sections with
support contained in $K$ and where $K$ runs through all compact
subsets of $M$.
 
So we declare the path components of $\Cal B= C^\infty(GL(TM,T^*M))$ for
the Whitney $C^\infty$-topology also to be open. We get a
topology which is finer than the Whitney topology, where each
connected component is homeomorphic to an open subset in
$\Cal B_c=C^\infty_c(L(TM,T^*M))$.
So $\Cal B=C^\infty(GL(TM,T^*M))$ is a smooth manifold
modeled on nuclear (LF)-spaces, and the tangent bundle is given
by $T\Cal B=\Cal B\x\Cal B_c$.
 
\subheading{\nmb.{1.2}. Remarks}
The main reference for the infinite dimensional manifold
structures is \cite{Michor, 1980}. But the
differential calculus used there is not completely up to date,
the reader should consult \cite{Fr\"olicher, Kriegl, 1988}, whose
calculus is more natural and much easier to apply.  There
a mapping between locally convex spaces is smooth if and only if
it maps smooth curves to smooth curves. See also
\cite{Kriegl, Michor, 1990} for a setting for real analytic
mappings along the same lines and applications to manifolds of
mappings.
 
As a final remark let us add that the differential structure on
the space $\Cal B$ of non degenerate bilinear structures is not completely
satisfying, if $M$ is not compact. In fact $C^\infty(L(TM,T^*M))$ is
a topological vector space with the compact $C^\infty$-topology,
but the space $\Cal B=C^\infty(GL(TM,T^*M))$ of non degenerate bilinear
structures is not open in it. Nevertheless, we will see later that the
exponential mapping for some pseudo riemannian metrics on $\Cal B$ is
defined also for some tangent vectors which are not in $\Cal B_c$.
This is an indication that the most natural setting for
manifolds of mappings is based on the compact
$C^\infty$-topology, but that one loses existence of charts. In
\cite{Michor, 1984} a setting for infinite dimensional manifolds
is presented which is based on an axiomatic structure of smooth
curves instead of charts.
 
\subheading{\nmb.{1.3}. The metrics} The tangent bundle of the
space 
$$\Cal B=C^\infty(GL(TM,T^*M))$$ 
of bilinear structures is
$$T\Cal B=\Cal B\x\Cal B_c=C^\infty(GL(TM,T^*M))\x\ 
C^\infty_c(L(TM,T^*M)).$$
Then $b\in\Cal B$ induces two fiberwise bilinear forms on
$L(TM,T^*M)$ which are given by $(h,k)\mapsto\tr(b\i hb\i k)$ and
$(h,k)\mapsto\tr(b\i h)\tr(b\i k)$.
We split each endomorphism $H=b\i h:TM\to TM$ into its trace free part
$H_0 := H - \frac{\tr(H)}{\dim M}Id$ and its trace part which
simplifies some formulas later on. Thus we have
$\tr(b\i hb\i k)= \tr((b\i h)_0(b\i k)_0) +
\frac1{\dim M}\tr(b\i h)\tr(b\i k)$. The structure $b$ also induces a
volume density on the base manifold $M$ by the local formula
$$\Vol(b)=\sqrt{|\det (b_{ij})}|\,|dx_1\wedge \dots\wedge dx_n|.$$
For each real $\al$ we have a smooth symmetric bilinear form
on $\Cal B$, given by
$$G^\al_b(h,k)=\int_M(\tr((b\i h)_0(b\i k)_0)
     +\al\tr(b\i h)\tr(b\i k))\Vol(b).$$
It is invariant under the action of the diffeomorphism
group $\Diff(M)$ on the space $\Cal B$ of bilinear structures. The integral is
defined since $h$, $k$ have compact support.
For $n=\dim M$ we have
$$G_b(h,k) := G^{1/n}_b(h,k) = \int_M\tr(b\i hb\i k)\Vol(b),$$
which for positive definite $b$ is the usual metric on the space of
all Riemannian metrics considered by \cite{Ebin, 1970},
\cite{Freed, Groisser, 1989}, and \cite{Gil-Medrano, Michor, 1990}.
We will see below in
\nmb!{1.4} that for $\al\ne0$ it is weakly non degenerate, i.e.
$G^\al_b$ defines a linear injective mapping
from the tangent space $T_b\Cal B=\Cal B_c=C^\infty_c(L(TM,T^*M))$
into its dual $C^\infty_c(L(TM,T^*M))'$, the space of distributional
densities with values in the dual bundle. This linear mapping is,
however, never
surjective. So we have a one parameter family of pseudo Riemannian
metrics on the infinite dimensional space $\Cal B$. The use of the
calculus of
\cite{Fr\"olicher, Kriegl, 1988} makes it completely obvious that it
is smooth in all appearing variables.
 
\proclaim{\nmb.{1.4}. Lemma} For $h,k\in T_b\Cal B$ we have
$$\align & G^\al _b(h,k)=G_b(h+\frac{\al n-1}n \tr(b^{-1}h)b,k),\\
&G_b(h,k)=G^\al _b(h-\frac{\al n-1}{\al n^2}\tr(b^{-1}h)b,k),
\text{ if }\al\ne 0,
\endalign$$
where $n=\dim M$. The pseudo Riemannian metric $G^\al$ is weakly non
degenerate for all $\al\ne0$.
\endproclaim

\demo{Proof}
The first equation is an obvious reformulation of the definition,
the second follows since
$h\mapsto h-\frac{\al n-1}{\al n^2}\tr(b^{-1}h)b$ is the
inverse of the transform $h\mapsto h+\frac{\al n-1}n\tr(b^{-1}h)b$.
Since $\tr(b_x\i h_x(b_x\i h_x)^{t,g})>0$ if
$h_x\ne0$, where ${\ell}^{t,g}$ is the transposed of a linear mapping
with respect to an arbitrary fixed Riemannian metric $g$, we have
$$G_b(h,b(b\i h)^{t,g})=\int_M\tr(b\i h(b\i h)^{t,g})\Vol(b) > 0$$
if $h\ne0$. So $G$ is weakly non degenerate, and by the second
equation $G^\al$ is weakly non degenerate for $\al\ne0$.
\qed\enddemo
 
\subheading{\nmb.{1.5}. Remark} Since $G^\al$ is only a weak
pseudo Riemannian metric, all objects which are only implicitly given
a priori lie in the Sobolev completions of the relevant spaces.
In particular this applies to the formula
$$\align 2G^\al (\xi,\nabla^\al_\et\ze)
=&\xi G^\al(\et,\ze) + \et G^\al (\ze,\xi) - \ze G^\al (\xi,\et)\\
&+ G^\al ([\xi,\et],\ze) +G^\al ([\et,\ze],\xi) - G^\al ([\ze,\xi],\et),
\endalign$$
which a priori gives only uniqueness but not existence of the
Levi Civita covariant derivative. But we refer to \cite{Gil-Medrano, 
Michor, 1990, 2.1} for a careful explanation of the role of covariant 
derivatives etc.
 
\proclaim{\nmb.{1.6}. Lemma} For $x\in M$ the pseudo metric
on $GL(T_xM,T^*_xM)$ given by
$$\ga^\al_{b_x}(h_x,k_x) := \tr((b_x\i h_x)_0 (b_x\i k_x)_0) +
\al\tr(b_x\i h_x)\tr(b_x\i k_x)$$
has signature (the number of negative eigenvalues)
$\frac{n(n-1)}2$ for $\al>0$
and has signature
$(\frac{n(n-1)}2+1)$ for $\al<0$.
\endproclaim
 
\demo{Proof}
In the framing $H=b_x\i h_x$ and $K=b_x\i k_x$ we have to determine
the signature of the symmetric bilinear form
$H,K\mapsto \tr(H_0K_0)+\al\tr(H)\tr(K)$. Since the signature is constant
on connected components we have to determine it only for $\al=\frac1n$ and
$\al=\frac1n-1$.
 
For $\al=\frac1n$ we note first that on the space of matrices
$H,K\mapsto \tr(HK^t)$ is positive
definite, and since the linear isomorphism $K\mapsto K^t$ has the
space of symmetric matrices as eigenspace for the eigenvalue 1, and
has the space of skew symmetric matrices as eigenspace for the
eigenvalue $-1$, we conclude that the signature is
$\frac{n(n-1)}2$ in this case.
 
For $\al=\frac1n-1$ we proceed as follows:
On the space of matrices with zeros on the main diagonal the
signature of $H,K\mapsto \tr(HK)$ is
$\frac{n(n-1)}2$ by the argument above and the form
$H,K\mapsto -\tr(H)\tr(K)$ vanishes. On the space
of diagonal matrices which we identify with $\Bbb R^n$
the whole bilinear form is given by
$\langle x,y\rangle = \sum_ix^iy^i-(\sum_ix^i)(\sum_iy^i)$.
Let $(e_i)$ denote the standard basis of $\Bbb R^n$ and put
$a_1:=\frac1n(e_1+\dots+e_n)$ and
$a_i:=\frac1{\sqrt{i-1+(i-1)^2}}(e_1+\dots+ e_{i-1}-(i-1)e_i)$ for $i>1$.
Then $\langle a_1,a_1\rangle = -1+\frac1n$ and for $i>1$ we get
$\langle a_i,a_j\rangle=\de_{i,j}$. So the signature there is
$1$.
\qed\enddemo

\heading \totoc\nmb0{2}. Geodesics, Levi Civita connection, and curvature
\endheading
 
\subheading{\nmb.{2.1}} Let $t\mapsto b(t)$ be a smooth curve in
$\Cal B$: so $b:\Bbb R\x M\to GL(TM,T^*M)$ is smooth and by the choice of
the topology on $\Cal B$
made in \nmb!{1.1} the curve
$b(t)$ varies only in a compact subset of $M$, locally
in $t$, by \cite{Michor, 1980, 4.4.4, 4.11, and 11.9}.
Then its energy is given by
$$\align E_a^b(b):&= \tfrac12\int_a^bG^\al_b(b_t,b_t)dt\\
&=\tfrac12\int_a^b\int_M
     \left(\tr((b\i b_t)_0 (b\i b_t)_0)+\al\tr(b\i b_t)^2\right)\Vol(b)\,dt,
\endalign$$
where $b_t=\frac{\partial}{\partial t}b(t)$.
 
Now we consider a variation of this curve, so we assume now that
$(t,s)\mapsto b(t,s)$ is smooth in all variables and locally in
$(t,s)$ it only varies within a compact subset in $M$ --- this is
again the
effect of the topology chosen in \nmb!{1.1}. Note that $b(t,0)$ is
the old $b(t)$ above.
 
\proclaim{\nmb.{2.2}. Lemma}
In the setting of
\nmb!{2.1} we
have the first variation formula
$$\align \dd s  \vert  _0 &E(G ^\al)^{a_1}_{a_0}(b(\quad,s))
    =G_b^\al (b_t,b_s) \vert_{t=a_0}^{t=a_1}+\\
&\quad +\int_{a_0}^{a_1} G ( -b_{tt}+b_tb^{-1}b_t
    +\frac 14 \tr(b^{-1}b_tb^{-1}b_t)b-\frac 12 \tr(b^{-1}b_t)b_t+\\
&\qquad \qquad  +\al \, (-\tr(b^{-1}b_{tt})-\frac 14 \tr(b^{-1}b_t)^2+
    \tr(b^{-1}b_tb^{-1}b_t))b,b_s ) \;dt=\\
&=G^\al_b(b_t,b_s)\vert_{t=a_0}^{t=a_1}+\\
&\quad +\int_{a_0}^{a_1}G^\al( -b_{tt} + b_tb^{-1}b_t-
    \frac 12 \tr(b^{-1}b_t)b_t+\frac1{4\al n}\tr(b^{-1}_tb^{-1}b_t)b+\\
&\qquad\qquad\qquad\qquad+\frac{\al n-1}{4\al n^2}\tr(b^{-1}b_t)^2b,b_s) \;dt
\endalign$$
\endproclaim
 
\demo{Proof}
We may interchange $\frac{\partial}{\partial s}|_0$ with the
first integral describing the energy in \nmb!{2.1}
since this is finite dimensional analysis, and we
may interchange it with the second one, since $\int_M$ is a
continuous linear functional on the space of all smooth
densities with compact support on $M$, by the chain rule.
Then we use that $\tr_*$ is linear and continuous,
$d(\Vol)(b)h= \frac12\tr(b\i h)\Vol(b)$, and that
$d((\quad)\i)_*(b)h= -b\i h b\i$ and partial integration.
\qed\enddemo
 
\subheading{\nmb.{2.3}. The geodesic equation}
By lemma
\nmb!{2.2} the curve $t\mapsto b(t)$ is a geodesic if
and only if we have
$$\align 
b_{tt} &=   b_tb^{-1}b_t-\frac 12 \tr(b^{-1}b_t)b_t
    +\frac1{4\al n}\tr(b^{-1}b_tb^{-1}b_t)b
    +\frac{\al n-1}{4\al n^2}\tr(b^{-1}b_t)^2b.\\
&= \Ga_b(b_t,b_t),
\endalign$$
where the {\it $G^\al$-Christoffel symbol}
$\Ga^\al \:\Cal B \x \Cal B_c \x \Cal B_c \to \Cal B_c$  is
given by symmetrisation
$$\align 
\Ga ^\al _b(h,k) &=\frac 12 hb^{-1}k+\frac 12 kb^{-1}h
     -\frac 14 \tr(b^{-1}h)k-\frac 14 \tr(b^{-1}k)h+\\
&\quad +\frac1{4\al n} \tr(b^{-1}hb^{-1}k)b
    +\frac{\al n-1}{4\al n^2}\tr(b^{-1}h)\tr(b^{-1}k)b.
\endalign$$
The sign of $\Ga^\al$ is chosen in such a way that the horizontal
subspace of $T^2\Cal B$ is parameterized by $(x,y;z,\Ga_x(y,z))$.
If instead of the obvious framing we use
$T\Cal B=\Cal B\x\Cal B_c \ni (b,h)\mapsto (b,b\i h)=:(b,H)\in
\{b\}\x C^\infty_c(L(TM,TM))$, the Christoffel symbol looks like
$$\align
\overline {\Ga^\al}_b(H,K)&=
     \frac 12 (HK+KH)-\frac 14 \tr(H)K-\frac 14 \tr(K)H\\
&\quad+\frac1{4\al n} \tr(HK)Id+\frac{\al n-1}{4\al n^2}\tr(H)\tr(K).
\endalign$$
and the $G^\al$-geodesic equation for $B(t):= b\i b_t$ becomes
$$ B_t = \dd t (b^{-1}b_t) = \frac1{4\al n}\tr(BB)Id
    - \frac 12 \tr(B)B+\frac{\al n-1}{4\al n^2}\tr(B)^2Id.$$
 
\subheading{\nmb.{2.4}. The curvature}
For vector fields $X$, $Y\in \X(\Cal N)$ and a vector
field $s:\Cal N\to T\M$ along $f:\Cal N\to \M$ we have
$$R(X,Y)s = (\nabla_{[X,Y]} - [\nabla_X,\nabla_Y])s 
    = (K\o TK - K\o TK\o \ka_{T\M})\o T^2s\o TX\o Y,$$
where $K:TTM\to M$ is the connector (see \cite{Gil-Medrano, Michor, 2.1})
and where the second formula
in local coordinates reduces to the usual formula
$$R(h,k)\ell = d\Ga(h)(k,\ell) - d\Ga(k)(h,\ell) -
\Ga(h,\Ga(k,\ell)) + \Ga(k,\Ga(h,\ell)).$$
A global derivation of this formula can be found in
\cite{Kainz, Michor, 1987}.
 
\proclaim{\nmb.{2.5}. Theorem} The curvature for the
pseudo Riemannian metric $G^\al$ on the manifold $\Cal B$ of all
non degenerate bilinear structures is given by
$$\align
b^{-1}R^\al_b(h,k)l&=\frac 14[[H,K],L]+\frac 1{16\al}(-\tr(HL)K+\tr(KL)H)+\\
&\quad+\frac{4\al n-3\al n^2+4n-4}{16\al n^2}(\tr(H)\tr(L)K-\\
&\quad-\tr(K)\tr(L)H)+\\
&\quad+\frac{4\al^2n^2-4\al n+\al n^2+3}{16\al n^2}
     (\tr(HL)\tr(K)Id\\
&\quad-\tr(KL)\tr(H)Id),
\endalign$$
where $H=b^{-1}h,K=b^{-1}k$ and $L=b^{-1}l$.
\endproclaim
\demo{Proof}
This is a long but direct computation.
\qed\enddemo
 
The geodesic equation can be solved explicitly and we have
 
\proclaim{\nmb.{2.6}. Theorem} Let $b^0\in\Cal B$ and
$h\in T_{b^0}\Cal B=\Cal B_c$.
Then the geodesic for the metric $G^\al$ in $\Cal B$ starting at $b^0$
in the direction of $h$ is the curve
$$\Exp^\al_{b^0}(th)=b^0 e^{(a(t)Id+b(t)H_0)},$$
where $H_0$ is the traceless part of $H:= (b^0)\i h$ (i.e.
$H_0=H-\frac{\tr(H)}n Id$) and where $a(t)=a_{\al,H}(t)$ and
$b(t)=b_{\al,H}(t)$ in $C^\infty(M)$ are defined as follows:
$$\align
a_{\al,H}(t)&=\frac 2n \log\left((1+\frac t4 \tr(H))^2
    +t^2\frac{\al\i}{16} \tr(H_0^2)\right)\\
b_{\al,H}(t)&=\cases\dsize
  \frac 4{\sqrt{\al\i\tr(H^2_0)}}
        \arctan\left(\frac{t\,\sqrt{\al\i\tr(H^2_0)}}{4+t\tr(H)}\right)
        &\text{for $\al\i\tr(H_0^2)> 0$}\\
  \frac 4{\sqrt{-\al\i\tr(H^2_0)}}
        \artanh\left(\frac{t\,\sqrt{-\al\i\tr(H^2_0)}}{4+t\tr(H)}\right)
        &\text{for $\al\i\tr(H_0^2)< 0$}\\
  \frac t{1+\frac t4 \tr(H)} &\text{for $\tr(H^2_0)=0$}
      \endcases
\endalign$$
Here $\arctan$ is taken to have values in $(-\frac \pi 2,\frac \pi 2)$
for the points of the basis manifold, where $\tr(H)\ge 0$, and on a point where
$\tr(H)< 0$ we define
$$\arctan\left(\frac {t\sqrt {\al\i\tr(H_0^2)}}{4+t\tr(H)}\right)=
\cases
\arctan \text { in $[0,\frac \pi 2)$}
                &\text {for $t\in [0,-\frac 4{\tr(H)})$}\\
\frac \pi 2     &\text {for $t=-\frac 4{\tr(H)}$}\\
\arctan \text { in $(\frac \pi 2 ,\pi)$}
                                   &\text {for $t\in
                                   (-\frac 4{\tr(H)},\infty)$}.
\endcases$$
To describe the domain of definition of the exponential mapping we
consider the sets
$$\align
Z^h\:&=\{x\in M :\frac1\al\tr_x(H_0^2)=0 \text{ and }\tr_x(H)<0\},\\
G^h\:&=\{x\in M :0>\frac1\al\tr_x(H_0^2)>-\tr_x(H)^2
     \text{ and }\tr_x(H)<0\}\\
&=\{x\in M:\al\ga(h,h)\lessgtr\ga^\al(h,h)\lessgtr0
     \text{ for }\al\lessgtr0, \tr_x(H)<0 \},\\
E^h\:&=\{x\in M :-\tr_x(H)^2=\frac1\al\tr_x(H_0^2)\text{ and }\tr_x(H)<0\}\\
&=\{x\in M :\ga^\al(h,h)=0 \text{ and }\tr_x(H)<0\},\\
L^h\:&=\{x\in M :-\tr_x(H)^2>\frac1\al\tr_x(H_0^2)\}\\
&=\{x\in M : \ga^\al(h,h)\gtrless0\text{ for }\al\lessgtr0 \},\\
\endalign$$
where $\ga(h,h)=\tr_x(H^2)$, and
$\ga^\al(h,h)=\tr_x(H_0^2)+\al\tr_x(H)^2$, see \nmb!{1.6},
are the integrands of $G_{b^0}(h,h)$ and $G^\al_{b^0}(h,h)$,
respectively.
Then we consider the numbers
$$\align
&z^h\:=\inf\{-\frac 4{\tr_x(H)}\:x\in Z^h\},\\
&g^h\:=\inf\{
4\frac {-\al\tr_x(H)-\sqrt{-\al\tr_x(H_0^2)}}{\tr_x(H_0^2)+\al\tr(H)^2}
\:x\in G^h\},\\
&e^h\:=\inf\{-\frac 2{\tr_x(H)}\:x\in E^h\},\\
&l^h\:=\inf\{
4\frac {-\al\tr_x(H)-\sqrt{-\al\tr_x(H_0^2)}}{\tr_x(H_0^2)+\al\tr(H)^2}
\:x\in L^h\},
\endalign$$
if the corresponding set is not empty, with value $\infty$ if the set
is empty.
Put $m^h\:=\inf \{z^h,g^h,e^h,l^h\}$.
Then $\Exp^\al_{b^0}(th)$ is maximally defined for $t\in[0,m^h)$.
\endproclaim

The second representations of the sets $G^h$, $L^h$, and $E^h$ 
clarifies how to take care of timelike, spacelike, and lightlike 
vector, respectively.
 
\demo{Proof}
The geodesical equation is very similar to that of the metric $G$ on 
the space of all Riemannian metrics, whose solution can be found in 
\cite{Freed, Groisser, 1989}, see also \cite{Gil-Medrano, Michor, 
1990}. The difference now is essentially that one should control the 
sign of various appearing constants. Here we use a slightly simpler 
method that unable us to deal only with scalar equations. Using 
$X(t):=g^{-1}g_t$ the geodesic equation reads as 
$$X'=-\frac 12 \tr(X)X+\frac1{4\al n} \tr(X^2)Id+ 
    \frac{\al n-1}{4\al n^2}\tr(X)^2Id,$$ 
and it is easy to see that a solution  $X$ satisfies
$$X'_0=-\frac 12 \tr(X)X_0.$$
Then $X(t)$ is in the plane generated by $H_0$ and $Id$ for 
all $t$ and the solution has 
the form $g(t)=b^0\exp(a(t)Id+b(t)H_0).$ 
Since $g_t=g(t)(a'(t)Id+b'(t)H_0)$ we have
$$\align X(t)&=a'(t)Id+b'(t)H_0 \text { and}
\\ X'(t)&=a''(t)Id+b''(t)H_0,
\endalign$$ 
and the geodesic equation becomes
$$\align
a''(t)Id+b''(t)H_0=&-\frac 12 na'(t)(a'(t)Id+b'(t)H_0)+\\
&+\frac1{4\al n}(na'(t)^2+b'(t)^2\tr(H_0^2))Id+\\
&+\frac{\al n-1}{4\al n^2}(n^2a'(t)^2)Id.\endalign$$
We may assume that $Id$ and $H_0$ are linearly
independent; if not $H_0=0$ and $b(t)=0$.
Hence the geodesic equation reduces to the differential equation
$$\left\{
{\aligned
a''&=-\frac n4
(a')^2+\frac {\tr(H_0^2)}{4\al n}(b')^2\\
b''&=-\frac n2 a'\,b'\endaligned}\right.$$
with initial conditions
$a(0)=b(0)=0$, $a'(0)=\frac{\tr(H)}n$, and
$b'(0)=1$.
 
If we take $p(t) = \exp (\frac n2 a)$ it is easy to see that then $p$ 
should be a solution of $p''' = 0$ and from the initial conditions
$$p(t) = 1+\frac t2 \tr(H)+ \frac{t^2}{16} (\tr(H)^2 +
\al\i\tr(H_0^2)).$$
Using that the second equation becomes $b' = p\i$, and then $b$ is 
obtained just by computing the integral. The solutions are defined in 
$[0,m^h)$ where $m^h$ is the infimum over the support of $h$ of the 
first positive root of the polynomial $p$, if it exists, and $\infty$ 
otherwise. The description of $m^h$ is now a technical fact.
\qed\enddemo
 
\subheading{\nmb.{2.7}. The exponential mapping}
For $b^0\in GL(T_xM,T_x^*M)$ and $H=(b^0)\i h$ let $C_{b^0}$ be the subset
of $L(T_xM,T_x^*M)$ given by the union of the sets (compare with 
$Z^h$, $G^h$, $E^h$, $L^h$ from \nmb!{2.6})
$$\gather
\{h:\tr(H_0^2)=0,\tr(H)\le-4\},\\
\{h: 0>\frac1\al\tr(H_0^2)>-\tr(H)^2,
     4\frac {-\al\tr(H)-\sqrt{-\al\tr(H_0^2)}}
     {\tr(H_0^2)+\al\tr(H)^2}\le1, \tr(H)<0\},\\
\{h: -\tr(H)^2=\frac1\al\tr(H_0^2),\tr(H)<-2\},\\
\{h: -\tr(H)^2>\frac1\al\tr(H_0^2),
     4\frac {-\al\tr(H)-\sqrt{-\al\tr(H_0^2)}}
     {\tr(H_0^2)+\al\tr(H)^2}\le1\}^{\text{closure}}
\endgather$$
which by some limit considerations coincides with the union of the
following two sets:
$$\multline
\{h: 0>\frac1\al\tr(H_0^2)>-\tr(H)^2,\\
\shoveright{     4\frac {-\al\tr(H)-\sqrt{-\al\tr(H_0^2)}}
     {\tr(H_0^2)+\al\tr(H)^2}\le1, \tr(H)<0\}^{\text{closure}},}\\
\{h: -\tr(H)^2>\frac1\al\tr(H_0^2),
     4\frac {-\al\tr(H)-\sqrt{-\al\tr(H_0^2)}}
     {\tr(H_0^2)+\al\tr(H)^2}\le1\}^{\text{closure}}.
\endmultline$$
So $C_{b^0}$ is closed.
We consider the open sets $U_{b^0} := L(T_xM,T_x^*M)\setminus C_{b^0}$,
$U'_{b^0}:=\{(b^0)\i h:h\in U_{b^0}\}\subset L(T_xM,T_xM)$,
and finally the open sub fiber bundles over $GL(TM,T^*M)$
$$\align
U:&= \bigcup\left\{\{b^0\}\x U_{b^0}: b^0\in GL(TM,T^*M)\right\}\subset\\
&\quad\subset     GL(TM,T^*M)\x_ML(TM,T^*M),\\
U':&= \bigcup\left\{\{b^0\}\x U'_{b^0}: 
b^0\in GL(TM,T^*M)\right\}\subset\\
&\quad\subset     GL(TM,T^*M)\x_ML(TM,TM).
\endalign$$
Then we consider the mapping $\Ph:U\to GL(TM,T^*M)$ which is given by
the following composition
$$\multline
U @>\sharp >> U'@>\ph >> GL(TM,T^*M)\x_ML(TM,TM) @>Id\x_M\exp >> \\
     @>Id\x_M\exp>> GL(TM,T^*M)\x_MGL(TM,TM) @>\flat >> GL(TM,T^*M),
\endmultline$$
where $\sharp(b^0,h):= (b^0, (b^0)\i h)$ is a fiber respecting
diffeomorphism,
$\flat(b^0,H):= b^0H$ is a diffeomorphism for fixed $b^0$, and
where the other two mappings will be discussed below.
 
The usual fiberwise exponential mapping
$$\exp:L(TM,TM)\to GL(TM,TM)$$
is a diffeomorphism near the zero section, on the ball of radius
$\pi$ centered at zero in a norm on the Lie algebra for which the Lie
bracket is sub multiplicative, for example. If we fix a symmetric
positive definite inner product $g$, then $\exp$ restricts to a 
global diffeomorphism from the linear subspace of $g$-symmetric 
endomorphisms onto the open subset of matrices which are positive 
definite with respect to $g$. If $g$ has signature this is no longer 
true since then $g$-symmetric matrices may have non real eigenvalues.
 
On the open set of all matrices whose eigenvalues $\la$ satisfy 
$|\Im\la|<\pi$, the exponential mapping is a diffeomorphism, see 
\cite{Varadarajan, 1977}.
 
The smooth mapping $\ph:U'\to GL(TM,T^*M)\x_ML(TM,TM)$ is given by
$\ph(b^0,H):= (b^0,a_{\al,H}(1)Id + b_{\al,H}(1)H_0)$ (see theorem
\nmb!{2.6}). It is a diffeomorphism onto its image with the following
inverse:
$$\ps(H):=
\cases
\aligned
&\frac 4n\left( e^{\frac
{\tr(H)}4}\cos\left(\frac{\sqrt{\al\i\tr(H_0^2)}}4\right ) -1\right)Id+\\
&\frac 4{\sqrt{\al\i\tr(H_0^2)}} e^{\frac {\tr(H)}4}\sin\left(\frac
{\sqrt{\al\i\tr(H_0^2)}}4\right) H_0
\endaligned
&\text{if $\tr(H_0^2)\ne 0$}\\
\frac 4n\left(e^{\frac {\tr(H)}4}-1\right)Id
&\text{otherwise,}
\endcases$$
where $\cos$ is considered as a complex function,
$\cos(iz)=i\cosh(z)$.
 
The mapping $(pr_1,\Ph): U\to GL(TM,T^*M)\x_M GL(TM,T^*M)$ is a
diffeomorphism on an open neighborhood of the zero section in $U$.
 
\proclaim{\nmb.{2.8}. Theorem}
In the setting of
\nmb!{2.7} the
exponential mapping $\Exp^\al_{b^0}$ for the metric $G^\al$
is a real analytic mapping defined on the open subset
$$\Cal U_{b^0}:= \{h\in C^\infty_c(L(TM,T^*M)):(b^0,h)(M)\subset U\}$$
and it is given by
$$\Exp_{b^0}(h) = \Ph\o(b^0,h).$$
The mapping
$(\pi_{\Cal B},\Exp): T\Cal B\to \Cal B\x \Cal B$ is a real analytic
diffeomorphism from an open neighborhood of the zero section in
$T\Cal B$
onto an open neighborhood of the diagonal in $\Cal B\x \Cal B$.
$\Cal U_{b_0}$ is the maximal domain of definition for the exponential
mapping.
\endproclaim
 
\demo{Proof} Most assertions are easy consequences of the
considerations above. For real analyticity of $\Exp$
the proof of \cite{Gil-Medrano, Michor, 1990, 3.4} applies which made 
use of deep results from \cite{Kriegl, Michor, 1990}
\qed\enddemo
 
\heading\totoc\nmb0{3}. Some submanifolds of $\Cal B$ \endheading
 
\subheading{\nmb.{3.1}. Submanifolds of pseudo Riemannian metrics}
We denote by $\Cal M^{q}$ the space of all pseudo Riemannian metrics 
on the manifold $M$ of signature (the dimension of a maximal negative 
definite subspace) $q$. It is an open set in a closed locally affine 
subspace of $\Cal B$  and thus a splitting submanifold of it with 
tangent bundle $T\Cal M^{q}=\Cal M^{q}\x C^\infty_c(S^2T^*M)$.
 
We consider a geodesic
$b(t)=b^0 e^{(a(t)Id+b(t)H_0)}$
for the metric $G^\al$ in $\Cal B$ starting at $b^0$
in the direction of $h$ as in \nmb!{2.6}. If $b^0\in\Cal M^{q}$
then $h\in T_{b^0}\Cal M^{q}$ if and only if
$H=(b^0)\i h\in L_{\text{sym}, b^0}(TM,TM)$ is symmetric with respect
to the pseudo Riemannian metric $b^0$.
But then $e^{(a(t)Id+b(t)H_0)}\in L_{\text{sym}, b^0}(TM,TM)$
for all $t$ in the domain of definition of the geodesic, so
$b(t)$ is a curve of pseudo Riemannian metrics and thus of the
same signature $q$ as $b^0$. Thus we have
 
\proclaim{\nmb.{3.2}. Theorem} For each $q\le n=\dim M$
the submanifold $\Cal M^{q}$ of pseudo Riemannian metrics of
signature ${q}$ on $M$ is a geodesically closed submanifold of
$(\Cal B, G^\al)$ for each $\al\ne0$.
\endproclaim

\subheading {Remark} The geodesics of $(\Cal M^{0},G^\al)$ have been 
studied, for $\al = \frac 1n$, in \cite{Freed, Groisser, 1989},
\cite{Gil-Medrano, Michor, 1990} and from \nmb!{3.2} and \nmb!{2.6} 
we see that they are completely analogous for every positive $\al$.
 
For fixed $x\in M$ there exists a family of homothetic pseudo metrics 
on the finite dimensional manifold $S^2_+T_x^*M$ whose geodesics are 
given by the evaluation of the geodesics of $(\Cal M^{0},G^\al)$ (see
\cite{Gil-Medrano, Michor, 1990} for more details). When $\al$ is 
negative, it is not difficult to see, from \nmb!{3.2} and \nmb!{2.6} 
again, that a geodesic of $(\Cal M^{0},G^\al)$ is defined for all $t$ 
if and only if the initial velocity $h$ satisfies $\ga^\al(h,h)\le 0$ 
and $\tr H > 0$ at each point of $M$ and then the same is true for 
all the above pseudo metrics on $S^2_+T_x^*M.$ These results appear 
already in \cite{DeWitt, 1967} for $n = 3.$
 
\subheading {\nmb.{3.3}. The local signature of $G^\al$}
Since $G^\al$ operates in infinite dimensional spaces, the usual
definition of signature is not applicable.
But for fixed $g\in\Cal M^{q}$ the
signature of
$$\ga^\al_{g_x}(h_x,k_x) = \tr((g_x^{-1}h_x)_0(g_x^{-1}k_x)_0)
     + \al \tr(g_x^{-1}k_x)\tr(g_x^{-1}k_x)$$
on $T_g (S^2_{q}T^*_xM) = S^2T^*_xM$
is independent of $x\in M$ and the special
choice of $g\in \Cal M^{q}$. We will call it the
{\it local signature} of $G^\al$.
 
\proclaim{\nmb.{3.4}. Lemma}
The signature of the quadratic form of
\nmb!{3.3} is
$$Q(\al,q)=q(q-n)+{\cases 0 &\text{ for }\al>0\\
                          1 &\text{ for }\al<0.\endcases}$$
\endproclaim
 
This result is due to \cite{Schmidt, 1989}.
 
\demo{Proof} Since the signature is constant
on connected components we have to determine it only for $\al=\frac1n$ and
$\al=\frac1n-1$. In a basis for $TM$ and its dual basis for $T^*M$ 
the bilinear form
$h\in S^2T_x^*M$ has a symmetric matrix. If the basis is 
orthonormal for $g$ we have (for $A^t=A$ and $C^t=C$)
$$H = g\i h = \pmatrix -Id_q & 0 \\ 0 & Id_{n-q}\endpmatrix
     \pmatrix A & B \\ B^t & C\endpmatrix
= \pmatrix -A & -B \\ B^t & C\endpmatrix,$$
which describes the typical matrix in the space
$L_{\text{sym},g}(T_xM,T_xM)$ of all
$H\in L(T_xM,T_xM)$ which are symmetric with respect to $g_x$.
 
Now we treat the case $\al=\frac1n$.
The standard inner product $\tr(HK^t)$ is positive definite on
$L_{\text{sym},g}(T_xM,T_xM)$ and the linear mapping $K\mapsto K^t$
has an eigenspace of dimension $q(n-q)$ for the eigenvalue $-1$ in
it, and a complementary eigenspace for the eigenvalue 1.
So $\tr(HK)$ has signature $q(n-q)$.
 
For the case $\al=\frac1n-1$ we again split the space
$L_{\text{sym},g}(T_xM,T_xM)$ into the subspace with 0 on the
main diagonal, where $\ga^\al_g(h,k)=\tr(HK)$ and
where $K\mapsto K^t$ has again an eigenspace of
dimension $q(n-q)$ for the eigenvalue $-1$, and the space of diagonal
matrices. There $\ga^\al_g$ has signature $1$ as determined in
the proof of \nmb!{1.6}.
\qed\enddemo
 
\subheading{\nmb.{3.5}. The submanifold of almost symplectic 
structures}\newline
A 2-form $\om\in\Om^2(M)=C^\infty(\La^2T^*M)$  can be non degenerate
only if $M$ is of even dimension $\dim M=n=2m$.
Then $\om$ is non degenerate if and only if
$\om\wedge \dots\wedge \om=\om^m$ is nowhere vanishing.
Usually this latter $2m$-form is regarded as the volume form
associated with $\om$,
but a short computation shows that we have
$$\Vol(\om)=\tfrac1{m!}|\om^m|.$$
This implies $m\ph\wedge \om^{m-1}=\frac12\tr(\om\i\ph)\om^m$.
 
\proclaim{\nmb.{3.6}. Theorem} The space $\Om^2_{\text{nd}}(M)$ of
non degenerate 2-forms is a splitting geodesically closed submanifold
of $(\Cal B, G^\al)$ for each $\al\ne0$.
\endproclaim
 
\demo{Proof}
We consider a geodesic
$b(t)=b^0 e^{(a(t)Id+b(t)H_0)}$
for the metric $G^\al$ in $\Cal B$ starting at $b^0$
in the direction of $h$ as in \nmb!{2.6}. If
$b^0=\om\in\Om^2_{\text{nd}}(M)$  then
$h\in\Om_c^2(M)$ if and only if  $H=\om\i h$
is symmetric with respect to $\om$, since we have
$\om(HX,Y)=\langle \om\om\i hX,Y\rangle = \langle hX, Y\rangle =
h(X,Y)=-h(Y,X)=-\om(HY,X)=\om(X,HY)$. At a point $x\in M$ we may
choose a Darboux frame $(e_i)$ such that $\om(X,Y)=Y^tJX$ where
$J=\pmatrix 0 & Id \\ -Id & 0\endpmatrix$. Then
$h$ is skew if and only if $JH$ is a skew symmetric matrix in the
Darboux frame, or $JH=H^tJ$. Since $(e^A)^t=e^{A^t}$ the matrix
$e^{a(t)Id+b(t)H_0}$ has then the same property, $b(t)$ is skew for
all $t$. So $\Om^2_{\text{nd}}(M)$ is a geodesically closed
submanifold.
\qed\enddemo
 
\proclaim{\nmb.{3.7}. Lemma} For a non degenerate 2-form $\om$ the
signature of the quadratic form
$\ph\mapsto \tr(\om\i\ph\om\i\ph)$
on $\La^2T^*_xM$ is
$m^2-m$ for $\al>0$ and $m^2-m+1$ for $\al<0$.
\endproclaim
 
\demo{Proof}
Use the method of \nmb!{1.6} and \nmb!{3.4}; the description of the
space of matrices can be read of the proof of \nmb!{3.6}.
\qed\enddemo
 
\subheading{\nmb.{3.8}. Symplectic structures}
The space $\operatorname{Symp}(M)$ of all symplectic structures is a
closed submanifold of $(\Cal B, G^\al)$. For a compact manifold $M$ it 
is splitting by the Hodge decomposition Theorem. For $\dim M=2$ we 
have $\operatorname{Symp}(M)=\Om^2_{\text{nd}}(M)$, so it is 
geodesically closed. But for $\dim M\ge 4$ the submanifold 
$\operatorname{Symp}(M)$ is not geodesically closed. For 
$\om\in\operatorname{Symp}(M)$ and $\ph$, $\ps\in
T_\om\operatorname{Symp}(M)$ the Christoffel form 
$\Ga^\al_\om(\ph,\ps)$ is not closed in general. The direct approach 
would need variational calculus with a partial differential equation 
as constraint. This will be treated in another paper. 
 
\heading\totoc\nmb0{4}. Splitting the manifold of metrics
\endheading
 
The developments in this sections were ignited by a question posed by
Maria Christina Abbati. The second author wants
to thank her for her question.
 
\subheading{\nmb.{4.1}. Almost product structures with given vertical
distribution} Let $M$ be a smooth finite dimensional manifold,
connected for simplicity's sake, and let $V$ be a distribution on
it. We will denote by $(TM,\pi,M)$ the tangent bundle,
 by $(V,\pi_V,M)$ the vector subbundle determined by $V$,
and by $(N=TM/V,\pi_N,M)$ the normal bundle. Let  $i:V\hookrightarrow
TM$ denote the embedding of $V$ and $p:TM\twoheadrightarrow
N$ the epimorphism onto the normal bundle.
 
Let us recall that an almost product structure on a manifold $M$ is a
(1,1)-tensor field $P$ (i.e. $P\in C^\infty(L(TM,TM))$) such that
$P^2 = Id$. It is evident that an almost product structure $P$ on
$M$ induces a decomposition of $TM$ of the form $TM = \ker(P - Id)
\oplus\ker(P + Id)$. These subbundles are called vertical and
horizontal and will be denoted by $V^P$, $H^P$ respectively . We also
have in a natural way two projectors $v^P = {1\over 2}(P + Id)$ and
$h^P = {1\over 2}(Id - P)$, the vertical (over $V^P$), and the
horizontal (over $H^P$)  projections. The almost product structure
$P$ also determines a
monomorphism $C_P:N\to TM$, called the horizontal lifting, given by
$C_P\o p = h^P$; it is an isomorphism onto $H^P$ inverse of $p|H^P$.
 
For a given distribution $V$ in $M$ we will denote by ${\Cal
P}_V(M)$ the space of all almost product structures with vertical
$V$ (i.e. such that $V = V^P$). So, giving an element of ${\Cal
P}_V(M)$ is equivalent to choosing a subbundle of $TM$ supplementary
of $V$, this subbundle is given then by $\ker(P + Id)$.
 
\proclaim{\nmb.{4.2}. Proposition}
The space ${\Cal P}_V(M)$ of almost product structures with vertical 
distribution $V$ is a real
analytic manifold with trivial tangent bundle whose fiber is
$\{\xi\in C^\infty_c(L(TM,TM)) : \im\xi\subset V\subset\ker\xi\}$
\endproclaim
 
\demo{Proof} We topologize $C^\infty(L(TM,TM))$ in such a way that it
becomes a topological locally affine space whose model vector
space is the space $$C^\infty_c(L(TM,TM))$$ of sections with
compact support. Then 
$${\Cal P}_V(M) = \{P\in C^\infty(L(TM,TM)): 
    P^2 = Id\ ,\ \ker(P - Id) = V\}$$
is a closed locally affine subspace of
$C^\infty(L(TM,TM))$, and thus a real analytic manifold.
 
The tangent space at $P$ is given by $$\{\xi\in C^\infty_c(L(TM,TM)): 
V\subset\ker\xi \ ,\ \xi P + P\xi = 0\}.$$
Now, for a $(1,1)$-tensor field $\xi$ the conditions 
$V\subset\ker\xi$ and $\xi P + P\xi = 0$
are equivalent to $h^P\xi = 0$ and $\xi v^P = 0$. The last couple of
conditions can be written only in term of $V$ as 
$\im\xi\subset V\subset\ker\xi$.
\qed\enddemo
 
\subheading {\nmb.{4.3}} For each metric $g$ on $M$ we have a
canonical choice of a complementary of $V$, just by taking the
orthogonal with respect to that metric, $V^{\bot,g}$, that defines
an almost product structure given by $P|V = Id$ and $P|V^{\bot,g}
= -Id$ (This structure is such that $(g,P)$ is an almost product
Riemannian structure, i.e. $g(P\cdot,P\cdot) = g(\cdot, \cdot)$ ).
$g$ also determines a metric on the bundle $(V,\pi_V,M)$ simply by
restriction and a metric on the normal bundle $(N,\pi_N,M)$ as the
restriction to $V^{\bot,g}$ via the isomorphism given by the horizontal
lifting.
 
Conversely, given an element $P$ of ${\Cal P}_V(M)$ and metrics
$g_1\in \M (N)$ and $g_2\in \M (V)$ a metric on $M$ can be
defined by $g(\cdot ,\cdot ) = g_1(p\cdot , p\cdot ) +
g_2(v^P\cdot ,v^P\cdot )$. It is easy to see that a bijection
is then established between $\M (M)$ and ${\Cal P}_V(M) \times\M
(N)\times\M (V)$.
 
\proclaim {\nmb.{4.4}. Proposition} There is a real analytic
diffeomorphism
$$\M (M) \cong \M(N)\times\M (V)\times{\Cal P}_V(M) $$
\endproclaim
 
\demo {Proof} In order to show that the above bijection is in fact a
real analytic diffeomorphism it will be convenient to write the maps
in the following way:
 
Let $\Ph$ be the map from $\M (M)$ to $\M
(N)\times\M (V)\times{\Cal P}_V(M) $ and let $\Pi_1$, $\Pi_2$, $\Pi_3$
the projections. We identify each metric $g$ with its associated
mapping $g:TM\to T^*M$, so that $\M(M)\subset C^\infty(L(TM,T^*M))$.
 
 We let $g_V$ denote the restriction of the metric $g$ to the subbundle
$V$, associated to it is the vector bundle isomorphism $g_V = i^* g
i:V\to V^*$, where $i:V\to TM$ is the injection and
$i^*:T^*M\to V^*$ is its adjoint. Then $\Pi_2\circ\Ph (g) = g_V$.
 
It is easy to see that the associated almost product structure
described above is given by
$\Pi_3\circ\Ph (g) = 2ig_V^{-1}i^* g - Id$.
 
Let us denote $C_g$ the horizontal lifting determined by
$\Pi_3\circ\Ph (g)$ (then, $C_g p = Id - ig_V^{-1}i^* g$) and
$C_g^*$ its adjoint. Then, $\Pi_1\circ\Ph (g) =C_g^* g C_g$.
 
Now, the inverse of $\Ph$ is given by $\Ps
(g_1,g_2,P) = p^*g_1p + (v^P)^*g_2v^P.$
 
Thus both $\Ph$ and $\Ps$ are the push forward of
sections by a fiber respecting smooth mapping which is fiberwise
quadratic, so it extends to a fiberwise holomorphic mapping in a
neighborhood between the complexifications of the affine bundles in
question. By the argument used in the proof of \cite{Gil-Medrano, 
Michor, 3.4} they are real analytic.
\qed\enddemo
 
\subheading {\nmb.{4.5}} We have seen that any distribution on $M$ 
induces a 
product structure on $\M$. We consider now the Riemannian
manifold $(\M , G)$, and we are going to see that there exists a metric on
$\M (N)\times\M (V)$ and a family of metrics on ${\Cal P}_V(M) $ such
that $G$ is what is usually called a product manifold with varying
metric on the fibers, although it is not the product metric.
 
To show that we will need some formulas which are obtained by
straightforward computations.
 
For each $(g_1,g_2)\in\M (N)\times\M (V)$ we have the immersion
$\Ps_{(g_1,g_2)}: {\Cal P}_V(M)\to \M$. The tangent map, at a point
$P\in{\Cal P}_V(M)$, 
$$T_P\Ps_{(g_1,g_2)}:\{\xi\in
     C^\infty_c(L(TM,TM)) : \im\xi\subset V\subset\ker\xi\}
    \to C^\infty_c(S^2T^*M)$$
is given by
$T_P\Ps_{(g_1,g_2)}(\xi) = {1\over 2}\{(v^P)^*g_2\xi + \xi^*g_2v^P\}$.
 
For each $P\in{\Cal P}_V(M)$ we have the immersion
$\Ps_P: \M (N)\times\M (V)\to\M$. The tangent map, at a point
$(g_1,g_2)\in\M (N)\times\M (V)$,
$$T_{(g_1,g_2)}\Ps_P:C^\infty_c(S^2N^*)\times C^\infty_c(S^2V^*)\to
C^\infty_c(S^2T^*M)$$
is given by $T_{(g_1,g_2)}\Ps_P(h_1, h_2) = p^*h_1p + (v^P)^*h_2v^P$.
 
For each $g\in\M$ the tangent map, at $g$, of the submersion
$\Pi_3\circ\Ph$ is given by 
$T_g(\Pi_3\circ\Ph)(h) = 2v^P   g^{-1}h  h^P$ 
and the tangent map, at $g$, of the submersion
$(\Pi_1,\Pi_2)\circ\Ph$ is given by 
$(T_g(\Pi_1,\Pi_2)\circ\Ph)(h)  = (C^*_g h C_g, i^* h i).$
 
\subheading {\nmb.{4.6}} If $P$ is an almost product structure on a
manifold $M$, each element $H\in C^\infty(L(TM,TM))$ can be written
as $H = H_1 + H_2$ where we have
$H_1 = v^P H v^P + h^P H h^P$ and
$H_2 = v^P H h^P + h^P H v^P$. That gives a decomposition of
$$
C^\infty_c(L(TM,TM)) = \D_1(P)\oplus\D_2(P)
$$
where 
$$\align
\D_1(P) &= \{H\in C^\infty_c(L(TM,TM)) : H(V^P)\subset V^P \text{ and  }
     H(H^P)\subset H^P\}, \\
\D_2(P) &= \{H\in C^\infty_c(L(TM,TM)) : H(V^P)\subset H^P \text{ and }
     H(H^P)\subset V^P\}.
\endalign$$
 
\subheading {\nmb.{4.7}} Let us assume now that a distribution $V$
has been fixed on $M$, for $g\in\M$ and $i=1,2$, let us denote
$\D_i(g) = \{h\in C^\infty_c(S^2T^*M)\ ;\
g^{-1}h\in\D_i((\Pi_3\circ\Ph)(g))\}$. It is clear that if $h\in
C^\infty_c(S^2T^*M)$ and if we take $H = g^{-1}h$ then $h = h_1 +
h_2$ with $h_i = gH_i$, $i = 1,2$. It is straightforward to see that if
$H$ is $g$-symmetric, then $H_1$ and $H_2$ are also $g$-symmetric.
We have in that way two complementary distributions in $\M$.
 
\proclaim {\nmb.{4.8}. Proposition} These distributions are mutually
orthogonal with respect to the metric $G$ on $\M$. They are both
integrable, more precisely  the leaves are the slices of the product.
\endproclaim
 
\demo {Proof} $\D_1$ and $\D_2$ are orthogonal to each other because for
any $h,k\in T_g(\M)$ we have $\tr (H_1K_2) = \tr (H_2K_1) = 0$ and
then, by the definition of $G$, 
$$G_g(h,k) = G_g(h_1,k_1) +
G_g(h_2,k_2).$$
 
The tangent space $\Ps_{(g_1,g_2)}({\Cal P}_V(M))$ at the point
$g = \Ps_{(g_1,g_2)}(P)$ is the kernel of the tangent mapping
$T_g((\Pi_1,\Pi_2)\circ\Ph)$ which, by \nmb!{4.5}, is
exactly $\D_2(g)$. Analogously, the slice $\Ps_P(\M (N)\times\M(V))$
has as tangent at $g$ the space $\ker T_g(\Pi_3\circ\Ph)$ which,
again by \nmb!{4.5}, is equal to  $\D_1(g).$
\qed\enddemo

\proclaim {\nmb.{4.9}. Proposition} The distribution $\D_1$ gives 
rise to a
totally geodesic foliation. A  non constant geodesic of $\M$ issuing from
$g$ in the direction of a vector in $\D_2(g)$ has the property that it
never meets again the leaf of $\D_2$ passing through $g$ and that its
tangent vector is never again in $\D_2$.
\endproclaim
 
\demo {Proof}  Let $g\in\M$  and let $(g_1, g_2, P) = \Ph(g)$. The
geodesic starting from $g$ in direction of $h\in T_g\M$ is given by
\cite{Gil-Medrano, Michor, 3.2} or \nmb!{2.6}:
$$g(t):= \Exp^G_g(th) = g e^{(a(t)Id+b(t)H_0)}.$$
For the first assertion, if $h\in\D_1(g)$ we have
$e^{(a(t)Id+b(t)H_0)}(H^P)\subset H^P$, and
$e^{(a(t)Id+b(t)H_0)}(V)\subset V$; since
$e^{(a(t)Id+b(t)H_0)}$ is non singular both inclusions are in fact
equalities and, consequently, $V^{\bot,g(t)} =
H^P$. So $g(t)\in \Ps_P(\M (N)\times\M(V))$  for all
$t$.
 
Let us suppose now that $h\in\D_2(g)$ then $\tr H = 0$ and
consequently $H_0 = H$. From \cite{Gil-Medrano, Michor, 3.2} or \nmb!{2.6}
it is easy to see that $a(t)\ne 0$ unless the geodesic is constant.
For $v_1,v_2\in V$ we have
$g(t)(v_1,v_2) = e^{a(t)}g(v_1,v_2)$ and then
$g(t)\in\Ps_{(g_1,g_2)}({\Cal P}_V(M))$ only for $t=0$.
 
 From \cite{Gil-Medrano, Michor, 4.5}, or \nmb!{2.6} (for $\al = \frac 1n)$ 
we have that 
$$P(t):= g(t)\i g'(t) = e^{-\frac12 na(t)}
     \left(\frac{4\tr(H)+nt\tr(H^2)}{4n}Id + H_0\right),$$
where $n=\dim M$.
Now, for $h\in\D_2(g)$, we have $g'(t)\in\D_2(g(t))$ if and only if the
coefficient of $Id$ is zero. If this happens for some $t\ne 0$ then
$\tr(H^2)=0$ which implies that $H=0$.
\qed\enddemo
 
\subheading{\nmb.{4.10}} For each $P\in {\Cal P}_V(M)$, one can define 
a metric on $\M (N)\times\M(V)$ by $\check G^P = \Ps_P^*G$  and for 
each $(g_1,g_2)\in\M (N)\times\M(V)$ a metric on ${\Cal P}_V(M)$ can 
be defined by $\hat G^{g_1,g_2} =\Ps_{(g_1,g_2)}^*G.$ The next 
propositions are devoted to the study of these
metrics.
 
\proclaim {\nmb.{4.11}. Proposition} All the
metrics $\check G^P$ on $\M (N)\times\M(V)$ are the same.
\endproclaim
\demo {Proof} From \nmb!{4.5}, for $(h_1,h_2)$, $(k_1,k_2)\in
T_{(g_1,g_2)}\M (N)\times\M(V)$ we have
$$\check G^P_{(g_1,g_2)}((h_1,h_2), (k_1,k_2)) = G_g(h,k)$$
where 
$g = \Ps_{(g_1,g_2)}(P) = p^*g_1p + (v^P)^* g_2 (v^P)$, $h = p^*h_1p +
    (v^P)^* h_2 (v^P)$ and  
$k = p^*k_1p + (v^P)^* k_2 (v^P).$
 
It is easy to see that $H = g\i h$ is given by 
$H = C_P H_1 p + i H_2 v^P$ with $H_1 = g_1\i h_1$, 
$H_2 = g_2\i h_2$ and then, $HK = C_P H_1K_1 p + i H_2K_2 v^P$.
 
Now, if we take, at each point, a base of $TM$ which is obtained from 
basis of $V$ and $N$ via the maps $i$ and $C_P$ we see that $\tr HK = 
\tr H_1K_1 + \tr H_2K_2$ and then
$$G_g(h,k) = \int_M \{\tr H_1K_1 + \tr H_2K_2\}\Vol(g).$$
The integrand does not depend on $P$ and $\Vol(g)$ is also independent 
of $P$ because if we have a curve $g(t)$ in 
$\Ps_{(g_1,g_2)}({\Cal P}_V(M))$ then $g'(t)\in {\Cal D}_2(g(t))$ by 
the proof of \nmb!{4.8},
and so $\tr(g(t)\i g'(t)) = 0$. From the expression of 
$(\Vol(g(t)))'$ (see \nmb!{2.2}) we conclude that
$\Vol(g(t))$ is constant.
\qed\enddemo

We will denote this metric on $\M (N)\times\M(V)$ by $\check G$.
 
\proclaim {\nmb.{4.12}. Proposition} For each $(g_1,g_2)\in\M (N)\times
\M(V)$ the metric $\hat G^{g_1,g_2}$ on ${\Cal P}_V(M)$ is flat. Its
exponential map, at each point, is then given by
$$\Exp^{\hat G^{g_1,g_2}}_P(\xi) = P + \xi.$$
\endproclaim
 
\demo {Proof} From \nmb!{4.5} we see that for $\xi$,
$\et\in T_P{\Cal P}_V(M)$ we have
$$\hat G^{g_1,g_2}_P(\xi,\et) = G_g(h,k)$$
where 
$g = \Ps_{(g_1,g_2)}(P) = p^*g_1p + (v^P)^* g_2 (v^P)$, $2h =(v^P)^* g_2
    \xi + \xi^* g_2 (v^P)$ 
and $2k = (v^P)^* g_2 \et + \et^* g_2 (v^P).$
 
Now, $g_2 = i^*gi$ and then $2h =(v^P)^* g \xi + \xi^* g (v^P).$ 
Having in mind the definition of $v^P$ and the facts that $P^*g = g 
P$ and that, for $\xi\in T_P{\Cal P}_V(M)$, $P\xi = \xi$ we have that 
$2h = g \xi + \xi^* g$ and then $2H = \xi + g\i\xi^*g$; analogously 
$2K = \et + g\i\et^*g.$ Consequently
$4HK = \xi\et +\xi g\i\et^*g + g\i\xi^*g\et + g\i\xi^*\et^*g =
\xi g\i\et^*g + g\i\xi^*g\et$, the last equality because $\xi\et = 
0.$
 
The distribution $V$ is contained in the kernel of the mapping 
$\xi^*g_2\et:TM\to T^*M$; the annihilator of $V$ contains the image 
of this mapping.
So there is a unique mapping 
$\widetilde{\xi^*g_2\et}:N\to N^*$ such that $p^* 
\widetilde{\xi^*g_2\et} p = \xi^*g_2\et.$
 
Then 
$g\i\xi^*g\et = g\i\xi^*g_2\et  = C_Pg_1\i\widetilde{\xi^*g_2\et} p$ 
and by an argument similar to that in \nmb!{4.11} 
$\tr (g\i\xi^*g_2\et) = \tr (g_1\i\widetilde{\xi^*g_2\et}).$ So, we 
conclude that
$$\hat G^{g_1,g_2}_P(\xi,\et) = \frac 14 \int_M \{\tr 
    (g_1\i\widetilde{\xi^*g_2\et}) +  
    \tr (g_1\i\widetilde{\et^*g_2\xi})\}\Vol(g)$$
which is independent of $P$ because $\Vol(g)$ does not depend on $P$
as we have shown in \nmb!{4.11}.
 
So, all the metrics are flat and then it is immediate that geodesics 
are just straight lines.
\qed\enddemo 

\subheading {Remark} For an element 
$\xi\in T_P{\Cal P}_V(M)$, $\xi^2 = 0$ and
then $e^\xi = Id + \xi$ and geodesics can also be written  in the form 
$P(t) = Pe^{t\xi}.$
 
\proclaim {\nmb.{4.13}. Proposition} In the submanifold 
$\Ps_{(g_1,g_2)}({\Cal P}_V(M))$ the geodesic starting at 
$g = \Ps_{(g_1,g_2)}(P)$ in the direction of
$h\in \Cal D_2(g)$ is given by
$$g(t) = g(Id + t H + t^2 H^2 h^P).$$
\endproclaim
 
\demo{Proof} The splitting submanifold $\Ps_{(g_1,g_2)}({\Cal P}_V(M))$
of $(\M, G)$ with the restricted metric is isometric to 
$({\Cal P}_V(M),\hat G^{g_1,g_2})$. The geodesic of the submanifold 
$\Ps_{(g_1,g_2)}({\Cal P}_V(M))$
starting at $g = \Ps_{(g_1,g_2)}(P)$ in the direction of 
$h\in \Cal D_2(g)$ is given by
$$g(t) = \Ps_{(g_1,g_2)}(P(t)) = p^*g_1p + (v^{P(t)})^* g_2 (v^{P(t)}),$$
where $P(t) = P + 2tv^pg\i h h^P$, by \nmb!{4.5} and \nmb!{4.12}.
Using \nmb!{4.1} and \nmb!{4.3} we have 
$g_1 = C^*_g g C_g$, $g_2 = i^* gi$, $C_g p = h^P$ and then
$$\align
g(t)  &= (h^P)^* g h^P + (v^{P(t)})^* g (v^{P(t)})\\
&= (h^P)^* g h^P + (v^P)^* g (v^P) + t^2(v^PHh^P)^* g (v^PHh^P)\\
&\qquad + t\{(v^PHh^P)^*g v^P + (v^P)^*g (v^PHh^P)\}\\
&= g(Id + t\{h^PHv^P + v^PHh^P\} + t^2 h^PHv^PHh^P),
\endalign$$
the last equality because $h^P$, $v^P$, $H$ are $g$-symmetric. Finally,
recalling that $h\in\Cal D _2(g)$ we have
$$g(t) = g + th + t^2 gH^2 h^P.\qed$$
\enddemo
 
\proclaim {\nmb.{4.14}. Theorem} The map
$(\Pi_1,\Pi_2)\circ\Ph :(\M,G)\to (\M (N)\times\M(V), \check G)$
is a Riemannian submersion. In
fact, $(\M, G)$ is a product manifold with varying metric on the
fibers.
\endproclaim
\demo {Proof} It follows by straightforward computation that
$$G_g = ((\Pi_1,\Pi_2)\circ\Ph)^*\check G_{(g_1,g_2)} +
(\Pi_3\circ\Ph)^*\hat G^{g_1,g_2}_P.\qed$$
\enddemo

\enddocument